\documentclass[11pt,leqno]{article}
\normalfont
\renewcommand{\baselinestretch}{1.2}

\usepackage{amsfonts}

\newcommand{\QED}{\raisebox{0.5mm}{\fbox{\rule{0mm}{1.5mm}\ }}}

\newcounter{myfn}[page]

\newcommand{\MainTheorem}{Theorem 1.1}
\newcommand{\MainCorollary}{Corollary 1.2}

\newcommand{\ChainProductVersion}{Theorem 2.1}
\newcommand{\SplittingBasicsProposition}{Proposition 2.2} 

\newcommand{\ComponentsAreChains}{Proposition 3.1}

\newcommand{\FundPosets}{Figure 4.1}
\newcommand{\FundLatticeIdealsFigure}{Figure 4.2}
\newcommand{\GridPosets}{Figure 4.3}
\newcommand{\GridPosetsII}{Figure 4.4}
\newcommand{\AllFigures}{Figures 4.1-4.4}
\newcommand{\SomeFigures}{Figures 4.1, 4.3, and 4.4}
\newcommand{\StructureMaxResult}{Proposition 4.1}

\newcommand{\GridPosetColoring}{Theorem 5.1}

\newcommand{\StartSection}{1}
\newcommand{\BackgroundSection}{2}
\newcommand{\TwoColorGridSection}{3}
\newcommand{\SemistandardSection}{4}
\newcommand{\ColorFunctionSection}{5}
\newcommand{\ProofSection}{6}

\newcommand{\selt}{\mathbf{s}} \newcommand{\telt}{\mathbf{t}}
 
 \newcommand{\xelt}{\mathbf{x}}
\newcommand{\yelt}{\mathbf{y}} 

\newcommand{\myA}{\mbox{\sffamily A}}
\newcommand{\mysmallA}{\mbox{\footnotesize \sffamily A}}

\newcommand{\myC}{\mbox{\sffamily C}}
\newcommand{\mysmallC}{\mbox{\footnotesize \sffamily C}}

\newcommand{\myG}{\mbox{\sffamily G}}
\newcommand{\mysmallG}{\mbox{\footnotesize \sffamily G}}

\newcommand{\maxelt}{\mbox{\bf max}}
\newcommand{\EuclideanE}{\mathfrak{E}}
\newcommand{\WGF}{\mbox{\sffamily WGF}}
\newcommand{\RGF}{\mbox{\sffamily RGF}}

\newcommand{\color}{\mathbf{color}}
\newcommand{\comp}{\mathbf{comp}}

\newcommand{\chain}{\mathbf{chain}}

\newcommand{\dichromatic}{two-color }

\newcommand{\Lba}{L_{\mathfrak{g}}^{\beta\alpha}(\lambda)}

\newcommand{\LAone}{L_{\mysmallA_{1} \oplus \mysmallA_{1}}}

\newcommand{\LAtwo}{L_{\mysmallA_{2}}}

\newcommand{\LBtwo}{L_{\mysmallC_{2}}}

\newcommand{\LGtwo}{L_{\mysmallG_{2}}}

\newcommand{\digriddelta}{(P, \leq_{_{P}}, 
\mathbf{chain}: P \longrightarrow [m], 
\color: P \longrightarrow \Delta)}

\newcommand{\myarrow}[1]{\stackrel{#1}{\rightarrow}}

\newcommand{\VertexForPosets}[2]{
\setlength{\unitlength}{1cm}
\begin{picture}(0,0)
\put(-0.38,0){
\begin{picture}(0,0)
\put(0,0){\circle*{0.15}} 
\put(-0.6,-0.1){\footnotesize $$}
\put(0.2,-0.1){\footnotesize $#2$}
\end{picture}
}
\end{picture}
}

\newcommand{\ZVertexForPosets}[2]{
\setlength{\unitlength}{1cm}
\begin{picture}(0,0)
\put(-0.38,0){
\begin{picture}(0,0)
\put(0,0){\circle*{0.15}} 
\put(-0.6,-0.1){\footnotesize $z_{#1}$}
\put(0.2,-0.1){\footnotesize $#2$}
\end{picture}
}
\end{picture}
}

\newcommand{\AOneAOneAlpha}{
\setlength{\unitlength}{1cm}
\begin{picture}(0.5,0.75)
\put(-0.45,0.2){\footnotesize $v_{1}$}
\put(0,0.25){\circle*{0.15}} 
\put(0.2,0.25){\footnotesize $\alpha$}
\end{picture}
}
\newcommand{\AOneAOneBeta}{
\setlength{\unitlength}{1cm}
\begin{picture}(0.5,0.75)
\put(-0.45,0.2){\footnotesize $v_{1}$}
\put(0,0.25){\circle*{0.15}} 
\put(0.2,0.25){\footnotesize $\beta$}
\end{picture}
}
\newcommand{\ATwoAlpha}{
\setlength{\unitlength}{1cm}
\begin{picture}(2,2)
\put(0,0.25){
\begin{picture}(2,2)
\put(0.5,-0.1){\footnotesize $v_{2}$}
\put(1,0){\circle*{0.15}} 
\put(1.2,-0.1){\footnotesize $\beta$}
\put(-0.5,0.9){\footnotesize $v_{1}$}
\put(0,1){\circle*{0.15}} 
\put(0.2,0.9){\footnotesize $\alpha$}
\put(1,0){\line(-1,1){1}} 
\end{picture}
}
\end{picture}
}
\newcommand{\ATwoBeta}{
\setlength{\unitlength}{1cm}
\begin{picture}(2,2)
\put(0,0.25){
\begin{picture}(2,2)
\put(0.5,-0.1){\footnotesize $v_{2}$}
\put(1,0){\circle*{0.15}} 
\put(1.2,-0.1){\footnotesize $\alpha$}
\put(-0.5,0.9){\footnotesize $v_{1}$}
\put(0,1){\circle*{0.15}} 
\put(0.2,0.9){\footnotesize $\beta$}
\put(1,0){\line(-1,1){1}} 
\end{picture}
}
\end{picture}
}
\newcommand{\BTwoAAlpha}{
\setlength{\unitlength}{1cm}
\begin{picture}(3,4)
\put(0,0.75){
\begin{picture}(3,3)
\put(1.5,-0.1){\footnotesize $v_{3}$}
\put(2,0){\circle*{0.15}} 
\put(2.2,-0.1){\footnotesize $\alpha$}
\put(0.5,0.9){\footnotesize $v_{2}$}
\put(1,1){\circle*{0.15}} 
\put(1.2,0.9){\footnotesize $\beta$}
\put(-0.5,1.9){\footnotesize $v_{1}$}
\put(0,2){\circle*{0.15}} 
\put(0.2,1.9){\footnotesize $\alpha$}
\put(2,0){\line(-1,1){2}} 
\end{picture}
}
\end{picture}
}
\newcommand{\BTwoBBeta}{
\setlength{\unitlength}{1cm}
\begin{picture}(4,4)
\put(1,0.25){
\begin{picture}(3,3.5)
\put(0.5,-0.1){\footnotesize $v_{4}$}
\put(1,0){\circle*{0.15}} 
\put(1.2,-0.1){\footnotesize $\beta$}
\put(-0.5,0.9){\footnotesize $v_{3}$}
\put(0,1){\circle*{0.15}} 
\put(0.2,0.9){\footnotesize $\alpha$}
\put(0.5,1.9){\footnotesize $v_{2}$}
\put(1,2){\circle*{0.15}} 
\put(1.2,1.9){\footnotesize $\alpha$}
\put(-0.5,2.9){\footnotesize $v_{1}$}
\put(0,3){\circle*{0.15}} 
\put(0.2,2.9){\footnotesize $\beta$}
\put(0,1){\line(1,1){1}} 
\put(1,0){\line(-1,1){1}} 
\put(1,2){\line(-1,1){1}} 
\end{picture}
}
\end{picture}
}
\newcommand{\GTwoAAlpha}{
\setlength{\unitlength}{1cm}
\begin{picture}(4,6)
\put(0.25,1.25){
\begin{picture}(4,6)
\put(2.5,-0.1){\footnotesize $v_{6}$}
\put(3,0){\circle*{0.15}} 
\put(3.2,-0.1){\footnotesize $\alpha$}
\put(1.5,0.9){\footnotesize $v_{5}$}
\put(2,1){\circle*{0.15}} 
\put(2.2,0.9){\footnotesize $\beta$}
\put(0.5,1.9){\footnotesize $v_{4}$}
\put(1,2){\circle*{0.15}} 
\put(1.2,1.9){\footnotesize $\alpha$}
\put(1.5,2.9){\footnotesize $v_{3}$}
\put(2,3){\circle*{0.15}} 
\put(2.2,2.9){\footnotesize $\alpha$}
\put(0.5,3.9){\footnotesize $v_{2}$}
\put(1,4){\circle*{0.15}} 
\put(1.2,3.9){\footnotesize $\beta$}
\put(-0.5,4.9){\footnotesize $v_{1}$}
\put(0,5){\circle*{0.15}} 
\put(0.2,4.9){\footnotesize $\alpha$}
\put(1,2){\line(1,1){1}} 
\put(3,0){\line(-1,1){2}} 
\put(2,3){\line(-1,1){2}} 
\end{picture}
}
\end{picture}
}
\newcommand{\GTwoBBeta}{
\setlength{\unitlength}{1cm}
\begin{picture}(4,8)
\put(0.25,0.25){
\begin{picture}(4,8)
\put(1.5,-0.1){\footnotesize $v_{10}$}
\put(2,0){\circle*{0.15}} 
\put(2.2,-0.15){\footnotesize $\beta$}
\put(0.5,0.9){\footnotesize $v_{9}$}
\put(1,1){\circle*{0.15}} 
\put(1.2,0.9){\footnotesize $\alpha$}
\put(1.5,1.9){\footnotesize $v_{8}$}
\put(2,2){\circle*{0.15}} 
\put(2.2,1.9){\footnotesize $\alpha$}
\put(0.5,2.9){\footnotesize $v_{6}$}
\put(1,3){\circle*{0.15}} 
\put(1.2,2.9){\footnotesize $\beta$}
\put(2.5,2.9){\footnotesize $v_{7}$}
\put(3,3){\circle*{0.15}} 
\put(3.2,2.9){\footnotesize $\alpha$}
\put(-0.5,3.9){\footnotesize $v_{4}$}
\put(0,4){\circle*{0.15}} 
\put(0.2,3.9){\footnotesize $\alpha$}
\put(1.5,3.9){\footnotesize $v_{5}$}
\put(2,4){\circle*{0.15}} 
\put(2.2,3.9){\footnotesize $\beta$}
\put(0.5,4.9){\footnotesize $v_{3}$}
\put(1,5){\circle*{0.15}} 
\put(1.2,4.9){\footnotesize $\alpha$}
\put(1.5,5.9){\footnotesize $v_{2}$}
\put(2,6){\circle*{0.15}} 
\put(2.2,5.9){\footnotesize $\alpha$}
\put(0.5,6.9){\footnotesize $v_{1}$}
\put(1,7){\circle*{0.15}} 
\put(1.2,6.9){\footnotesize $\beta$}
\put(1,1){\line(1,1){2}} 
\put(1,3){\line(1,1){1}} 
\put(0,4){\line(1,1){2}} 
\put(2,0){\line(-1,1){1}} 
\put(2,2){\line(-1,1){2}} 
\put(3,3){\line(-1,1){2}} 
\put(2,6){\line(-1,1){1}} 
\end{picture} 
}
\end{picture}
}

\newcommand{\AoneAlphaIdeals}{
\setlength{\unitlength}{1.1cm}
\begin{picture}(3,3)
\put(0.25,2.5){$\LAone(1,0)$}
\put(1.5,1){\circle*{0.125}}
\put(1.5,2){\circle*{0.125}}
\put(1.7,0.95){\footnotesize $\emptyset$}
\put(1.6,1.95){\footnotesize $\langle 1 \rangle$}
\put(1.5,1){\line(0,1){1}}
\put(1.45,1.45){\footnotesize $\alpha$}
\end{picture}
}

\newcommand{\AoneBetaIdeals}{
\setlength{\unitlength}{1.1cm}
\begin{picture}(3,3)
\put(0.25,2.5){$\LAone(0,1)$}
\put(1.5,1){\circle*{0.125}}
\put(1.5,2){\circle*{0.125}}
\put(1.7,0.95){\footnotesize $\emptyset$}
\put(1.6,1.95){\footnotesize $\langle 1 \rangle$}
\put(1.5,1){\line(0,1){1}}
\put(1.45,1.45){\footnotesize $\beta$}
\end{picture}
}

\newcommand{\AtwoAlphaIdeals}{
\setlength{\unitlength}{1.1cm}
\begin{picture}(3,4)
\put(0.25,3.5){$\LAtwo(1,0)$}
\put(1.5,1){\circle*{0.125}}
\put(1.5,2){\circle*{0.125}}
\put(1.5,3){\circle*{0.125}}
\put(1.7,0.95){\footnotesize $\emptyset$}
\put(1.6,1.95){\footnotesize $\langle 2 \rangle$}
\put(1.6,2.95){\footnotesize $\langle 1 \rangle$}
\put(1.5,1){\line(0,1){2}}
\put(1.45,1.45){\footnotesize $\beta$}
\put(1.45,2.45){\footnotesize $\alpha$}
\end{picture}
}

\newcommand{\AtwoBetaIdeals}{
\setlength{\unitlength}{1.1cm}
\begin{picture}(3,4)
\put(0.25,3.5){$\LAtwo(0,1)$}
\put(1.5,1){\circle*{0.125}}
\put(1.5,2){\circle*{0.125}}
\put(1.5,3){\circle*{0.125}}
\put(1.7,0.95){\footnotesize $\emptyset$}
\put(1.6,1.95){\footnotesize $\langle 2 \rangle$}
\put(1.6,2.95){\footnotesize $\langle 1 \rangle$}
\put(1.5,1){\line(0,1){2}}
\put(1.45,1.45){\footnotesize $\alpha$}
\put(1.45,2.45){\footnotesize $\beta$}
\end{picture}
}

\newcommand{\BtwoAlphaIdeals}{
\setlength{\unitlength}{1.1cm}
\begin{picture}(3,5)
\put(0.25,4.5){$\LBtwo(1,0)$}
\put(1.5,1){\circle*{0.125}}
\put(1.5,2){\circle*{0.125}}
\put(1.5,3){\circle*{0.125}}
\put(1.5,4){\circle*{0.125}}
\put(1.7,0.95){\footnotesize $\emptyset$}
\put(1.6,1.95){\footnotesize $\langle 3 \rangle$}
\put(1.6,2.95){\footnotesize $\langle 2 \rangle$}
\put(1.6,3.95){\footnotesize $\langle 1 \rangle$}
\put(1.5,1){\line(0,1){3}}
\put(1.45,1.45){\footnotesize $\alpha$}
\put(1.45,2.45){\footnotesize $\beta$}
\put(1.45,3.45){\footnotesize $\alpha$}
\end{picture}
}

\newcommand{\BtwoBetaIdeals}{
\setlength{\unitlength}{1.1cm}
\begin{picture}(3,6)
\put(0.25,5.5){$\LBtwo(0,1)$}
\put(1.5,1){\circle*{0.125}}
\put(1.5,2){\circle*{0.125}}
\put(1.5,3){\circle*{0.125}}
\put(1.5,4){\circle*{0.125}}
\put(1.5,5){\circle*{0.125}}
\put(1.7,0.95){\footnotesize $\emptyset$}
\put(1.6,1.95){\footnotesize $\langle 4 \rangle$}
\put(1.6,2.95){\footnotesize $\langle 3 \rangle$}
\put(1.6,3.95){\footnotesize $\langle 2 \rangle$}
\put(1.6,4.95){\footnotesize $\langle 1 \rangle$}
\put(1.5,1){\line(0,1){4}}
\put(1.45,1.45){\footnotesize $\beta$}
\put(1.45,2.45){\footnotesize $\alpha$}
\put(1.45,3.45){\footnotesize $\alpha$}
\put(1.45,4.45){\footnotesize $\beta$}
\end{picture}
}

\newcommand{\GtwoAlphaIdeals}{
\setlength{\unitlength}{1.1cm}
\begin{picture}(3,8)
\put(0.25,7.5){$\LGtwo(1,0)$}
\put(1.5,1){\circle*{0.125}}
\put(1.5,2){\circle*{0.125}}
\put(1.5,3){\circle*{0.125}}
\put(1.5,4){\circle*{0.125}}
\put(1.5,5){\circle*{0.125}}
\put(1.5,6){\circle*{0.125}}
\put(1.5,7){\circle*{0.125}}
\put(1.7,0.95){\footnotesize $\emptyset$}
\put(1.6,1.95){\footnotesize $\langle 6 \rangle$}
\put(1.6,2.95){\footnotesize $\langle 5 \rangle$}
\put(1.6,3.95){\footnotesize $\langle 4 \rangle$}
\put(1.6,4.95){\footnotesize $\langle 3 \rangle$}
\put(1.6,5.95){\footnotesize $\langle 2 \rangle$}
\put(1.6,6.95){\footnotesize $\langle 1 \rangle$}
\put(1.5,1){\line(0,1){6}}
\put(1.45,1.45){\footnotesize $\alpha$}
\put(1.45,2.45){\footnotesize $\beta$}
\put(1.45,3.45){\footnotesize $\alpha$}
\put(1.45,4.45){\footnotesize $\alpha$}
\put(1.45,5.45){\footnotesize $\beta$}
\put(1.45,6.45){\footnotesize $\alpha$}
\end{picture}
}

\newcommand{\GtwoBetaIdeals}{
\setlength{\unitlength}{1.1cm}
\begin{picture}(5,12)
\put(1.75,11.5){$\LGtwo(0,1)$}
\put(5,6){\circle*{0.125}}
\put(2,5){\circle*{0.125}}
\put(2,7){\circle*{0.125}}
\put(3,1){\circle*{0.125}}
\put(3,2){\circle*{0.125}}
\put(3,3){\circle*{0.125}}
\put(3,4){\circle*{0.125}}
\put(3,6){\circle*{0.125}}
\put(3,8){\circle*{0.125}}
\put(3,9){\circle*{0.125}}
\put(3,10){\circle*{0.125}}
\put(3,11){\circle*{0.125}}
\put(4,5){\circle*{0.125}}
\put(4,7){\circle*{0.125}}
\put(3.15,10.95){\footnotesize $\langle 1 \rangle$}
\put(3.15,9.95){\footnotesize $\langle 2 \rangle$}
\put(3.15,8.95){\footnotesize $\langle 3 \rangle$}
\put(3.15,7.95){\footnotesize $\langle 4,5 \rangle$}
\put(4.15,7){\footnotesize $\langle 4,7 \rangle$}
\put(1.45,6.95){\footnotesize $\langle 5 \rangle$}
\put(2.15,5.95){\footnotesize $\langle 6,7 \rangle$}
\put(5.1,5.95){\footnotesize $\langle 4 \rangle$}
\put(1.45,4.95){\footnotesize $\langle 7 \rangle$}
\put(4.15,4.9){\footnotesize $\langle 6 \rangle$}
\put(3.15,3.95){\footnotesize $\langle 8 \rangle$}
\put(3.1,2.95){\footnotesize $\langle 9 \rangle$}
\put(3.1,1.95){\footnotesize $\langle 10 \rangle$}
\put(3.2,0.95){\footnotesize $\emptyset$}
\put(3,1){\line(0,1){3}}
\put(3,8){\line(0,1){3}}
\put(3,4){\line(-1,1){1}}
\put(3,4){\line(1,1){1}}
\put(2,5){\line(1,1){2}}
\put(4,5){\line(-1,1){2}}
\put(2,7){\line(1,1){1}}
\put(4,7){\line(-1,1){1}}
\put(4,5){\line(1,1){1}}
\put(5,6){\line(-1,1){1}}
\put(2.95,10.45){\footnotesize $\beta$}
\put(2.95,9.45){\footnotesize $\alpha$}
\put(2.95,8.45){\footnotesize $\alpha$}
\put(2.45,7.45){\footnotesize $\alpha$}
\put(3.45,7.45){\footnotesize $\beta$}
\put(2.45,6.45){\footnotesize $\beta$}
\put(3.45,6.45){\footnotesize $\alpha$}
\put(4.45,6.45){\footnotesize $\alpha$}
\put(2.45,5.45){\footnotesize $\beta$}
\put(3.45,5.45){\footnotesize $\alpha$}
\put(4.45,5.45){\footnotesize $\alpha$}
\put(2.45,4.45){\footnotesize $\alpha$}
\put(3.45,4.45){\footnotesize $\beta$}
\put(2.95,3.45){\footnotesize $\alpha$}
\put(2.95,2.45){\footnotesize $\alpha$}
\put(2.95,1.45){\footnotesize $\beta$}
\end{picture}
}

\begin{document}

\newpage
\setcounter{page}{1} 
\renewcommand{\baselinestretch}{1}

\vspace*{-0.7in}
\hfill {\footnotesize This version: May 15, 2022}

\vspace*{-0.1in}
\hfill {\footnotesize First version: January 1, 2019}

\begin{center}
{\Large \bf Analogs of Schur functions for rank two Weyl groups\\ obtained from grid-like posets\\ }

L.\ Wyatt Alverson II, Robert G.\ Donnelly, Scott J.\ Lewis, and Robert Pervine

\vspace*{-0.05in} 
Department of Mathematics and Statistics, Murray State
University, Murray, KY 42071
\end{center}

\begin{abstract}
In prior work, the authors, along with M.\ McClard, R.\ A.\ Proctor, and N.\ J.\ Wildberger, studied certain distributive lattice models for the `Weyl bialternants' (aka `Weyl characters') associated with the rank two root systems$\,$/$\,$Weyl groups. 
These distributive lattices were uniformly described as lattices of order ideals taken from certain grid-like posets, although the arguments connecting the lattices to Weyl bialternants were case-by-case depending on the type of the rank two root system. 
Using this connection with Weyl bialternants, these lattices were shown to be rank symmetric and rank unimodal, and their rank generating functions were shown to have beautiful quotient-of-products expressions. 
Here, these results are re-derived from scratch using completely uniform and elementary combinatorial reasoning in conjunction with some combinatorial methodology developed elsewhere by the second-listed author. 

\begin{center}
{\small \bf Mathematics Subject Classification:}\ {\small 05E05 
(05A15, 05E10, 20F55, 17B10)}\\
{\small \bf Keywords:}\  {\footnotesize distributive lattice, poset of join irreducibles, rank generating function, rank two root system, Weyl group, Weyl bialternant, splitting poset, semisimple Lie algebra representation}
\end{center} 
\end{abstract}

{\bf \S \StartSection\ \ Introduction.} At first glance, the main result of \cite{ADLMPPW} appears to be Lie theoretic in nature: 
that paper shows how certain polynomial-like representation theory invariants can be realized as weight generating functions for certain specified families of beautiful diamond-colored distributive lattices. 
Each such invariant is the `Weyl bialternant' $\chi_{_{\lambda}}$ (sometimes also called a `Weyl character') associated with a highest weight $\lambda$ irreducible representation of some rank two semisimple Lie algebra $\mathfrak{g}$. 
The proof of the main result of \cite{ADLMPPW} made use of a close connection the authors found between the objects used to build the distributive lattices 
and some type-dependent tableaux developed by P.\ Littelmann \cite{Lit} in his study of generalizations of the Littlewood--Richardson rule for decomposing products of irreducible $\mathfrak{g}$-modules. 

Still, assuming some familiarity with root systems and their associated Weyl groups, the main result of \cite{ADLMPPW}, and its corollaries, can be stated and understood in purely combinatorial terms. 
Our re-statement of this result appears as \MainCorollary\ below. 
In this paper we re-derive this result as a consequence of a more general statement (\MainTheorem) whose proof is combinatorial and completely uniform, besides one part of the statement which is obviously particularized to the respective rank two root systems. 
Before we state this main result and its corollary, we provide some qualitative comments on the main ideas and structures involved. 

A Weyl bialternant is a quotient of alternating sums, each of which is skew-invariant under the action of a Weyl group associated to a finite root system. 
For example, when the Weyl group is the symmetric group, then this bialternant can be viewed as the quotient-of-determinants description of a classical Schur function.  
Weyl bialternants can be specialized in various ways. 
A certain $q$-specialization yields a symmetric and unimodal polynomial, expressed as a quotient of products, whose coefficients are positive integers which we will view as counts of combinatorial objects. 
So, letting $q \rightarrow 1$ yields a quotient-of-products expression for a related count of combinatorial objects. 

Here, as in \cite{ADLMPPW}, we want to use diamond-colored distributive lattices,  
realized as lattices of order ideals taken from certain grid-like posets (`two-color grid posets'), as models for Weyl bialternants associated to the rank two root systems$\,$/$\,$Weyl groups. 
In particular, the Weyl bialternants will be weight-generating functions for these lattices, and the aforementioned $q$-specializations will be the rank-generating functions.  

Certain families of such diamond-colored distributive lattices were presented in \cite{ADLMPPW}; these are called `semistandard lattices'. 
We regard the aforementioned rank-generating function result for semistandard lattices to be one of the combinatorial/enumerative highlights of our use of distributive lattices to model Weyl bialternants in rank two:  
In \cite{ADLMPPW}, it is noted how semistandard lattices provide an answer, for the rank two root systems, to a  question posed by R.\ P.\ Stanley in 1979 (see \cite{StanUnim}). 
For some other aspects and consequences of this overall approach, see \cite{ADLP} and particularly the introduction to \cite{ADLMPPW}. 

The rank two root systems are here denoted $\myA_{1} \oplus \myA_{1}$, $\myA_{2}$, $\myC_{2}$, and $\myG_{2}$. 
We use $\{\alpha,\beta\}$ as a basis of simple roots, with $\alpha$ short. 
The respective `Cartan matrices' $M_{\Phi}$ (cf.\ \S \BackgroundSection) depicted below encapsulate some distinguishing data for each rank two root system $\Phi$. 
\begin{center}
\begin{tabular}{|c|c|c|c|}
\hline 
$\myA_{1} \oplus \myA_{1}$ & $\myA_{2}$ & $\myC_{2}$ & $\myG_{2}$\\
\hline 
\hline
\rule[-7mm]{0mm}{16mm} $\displaystyle \left(\begin{array}{cc} 2 & 0 \\ 
0 & 2\end{array}\right)$ 
&  
\rule[-7mm]{0mm}{16mm} $\displaystyle \left(\begin{array}{cc} 2 & -1 \\ 
-1 & 2\end{array}\right)$ 
& 
\rule[-7mm]{0mm}{16mm} $\displaystyle \left(\begin{array}{cc} 2 & -1 \\ 
-2 & 2\end{array}\right)$ 
& 
\rule[-7mm]{0mm}{16mm} $\displaystyle \left(\begin{array}{cc} 2 & -1 \\ 
-3 & 2\end{array}\right)$\\
\hline 
\end{tabular}
\end{center}

In this paragraph we provide some descriptive comments about the semistandard lattices for each rank two root system. 
Each $\myA_{1} \oplus \myA_{1}$-semistandard lattice is just the product of two chains. 
The $\myA_{2}$-semistandard lattices are the Gelfand--Tsetlin lattices that arise in connection with the Gelfand--Tseltin bases for the irreducible representations of the semisimple Lie algebra $\mathfrak{sl}(3,\mathbb{C})$ associated to this root system; these can be easily described as natural orderings of certain `semistandard'  (or `column strict') tableaux, see e.g.\ \cite{HL} and references therein.  
Stanley \cite{StanUnim} and R.\ A.\ Proctor \cite{PrGZ} were among the first to study these lattices in connection with root systems, Weyl groups, and Lie algebra representations. 
The $\myC_{2}$- and $\myG_{2}$-semistandard lattices were discovered by the second-listed author as natural analogs of the Gelfand--Tsetlin lattices, and were first formally studied in master's theses by L.\ W.\ Alverson \cite{Alv} and M.\ McClard \cite{McClard} respectively. 
Initially, these semistandard lattices were obtained as natural orders on type-specific semistandard tableaux from \cite{KN} and \cite{Lit}. 
See \cite{ADLP} and \cite{ADLMPPW} for a detailed combinatorial description of all of these rank two semistandard tableaux. 

The idea of describing semistandard lattices uniformly across type in terms of their underlying posets of join irreducibles (`semistandard posets') was originated by Donnelly and N.\ J.\ Wildberger as part of a larger program that is still in development. 
It was discerned by the second-listed author that the semistandard tableaux associated with rank two root systems could be more uniformly described within the more general and fairly pleasant combinatorial environment of what we call `two-color grid posets' and that properties of two-color grid posets could help unify some of the initially case-oriented results obtained by the authors of \cite{ADLP} and \cite{ADLMPPW}. 
As in those papers, here we will gain advantage by viewing semistandard posets as particular instances of two-color grid posets.  
See \AllFigures\ for a depiction of some semistandard posets$\,$/$\,$lattices. 

With this background in mind, what follows is a statement of the main result of this paper. 
This statement uses some notation and terminology that are developed in the next two sections.  
But notice that all of the conclusions of the theorem statement follow from only two key combinatorial hypotheses: ({\sl i}) That the originating two-color grid poset has the `max property' and ({\sl ii}) that its associated diamond-colored distributive lattice is `$M_{\Phi}$-structured'. 
Besides equation (1), the statement and proof (see \S \ProofSection) nowhere depend upon the classification of rank two root systems. 

\noindent
{\bf \MainTheorem}\ \ {\sl Let $\Phi$ be a rank two root system with simple root basis $\{\alpha,\beta\}$ ($\alpha$ is short), associated fundamental weights $\{\omega_{\alpha},\omega_{\beta}\}$, and associated Cartan matrix $M_{\Phi}$. 
Let $P$ be a two-color grid poset with vertex color set $\{\alpha,\beta\}$. 
Let $L := \mathbf{J}_{color}(P)$, and let $\lambda = a\omega_{\alpha}+b\omega_{\beta}$ be the (dominant) weight of the unique maximal element of $L$.  
Suppose that $P$ has the max property and that $L$ is $M_{\Phi}$-structured. 
Then $L$ is a splitting distributive lattice for the $\Phi$-Weyl bialternant $\chi_{_{\lambda}}$. 
In particular} $\WGF(L) = \chi_{_{\lambda}}${\sl , $L$ is rank symmetric and rank unimodal, and} $\displaystyle \RGF(L,q) = \mbox{\rule[-6mm]{0mm}{15mm}\large $\displaystyle \prod_{\alpha \in \Phi^{+}}$}\frac{1-q^{\langle \lambda+\varrho,\alpha^{\vee} \rangle}}{1-q^{\langle \varrho,\alpha^{\vee} \rangle}}$. 
{\sl The latter identity particularizes as follows for the four distinct rank two root systems:}
\begin{equation}\RGF(L,q) = \left\{\begin{array}{cl}
\rule[-4mm]{0mm}{10mm}\frac{(1-q^{a+1})(1-q^{b+1})}{(1-q)(1-q)} & \mbox{ \sl if }\, \Phi = \myA_{1} \oplus \myA_{1}\\
\rule[-4mm]{0mm}{10mm}\frac{(1-q^{a+1})(1-q^{b+1})(1-q^{a+b+2})}{(1-q)(1-q)(1-q^{2})} & \mbox{ \sl if }\, \Phi = \myA_{2}\\
\rule[-4mm]{0mm}{10mm}\frac{(1-q^{a+1})(1-q^{b+1})(1-q^{a+b+2})(1-q^{a+2b+3})}
{(1-q)(1-q)(1-q^{2})(1-q^{3})} & \mbox{ \sl if }\, \Phi = \myC_{2}\\
\rule[-4mm]{0mm}{10mm}\frac{(1-q^{a+1})(1-q^{b+1})(1-q^{a+b+2})(1-q^{a+2b+3})
(1-q^{a+3b+4})(1-q^{2a+3b+5})}{(1-q)(1-q)(1-q^{2})(1-q^{3})
(1-q^{4})(1-q^{5})} & \mbox{ \sl if }\, \Phi = \myG_{2}
\end{array}\right.\end{equation} 

The next corollary follows easily from \MainTheorem\ (see \S \ProofSection\ below) and re-capitulates the main results of \cite{ADLMPPW}. 

\noindent 
{\bf \MainCorollary}\ \ {\sl In the notation of \MainTheorem, suppose now that $P$ is a $\Phi$-semistandard poset with associated $\Phi$-semistandard lattice $L = \mathbf{J}_{color}(P)$ such that the unique maximal element of $L$ has weight $\lambda  = a\omega_{\alpha}+b\omega_{\beta}$. 
Then $L$ is a splitting distributive lattice for the $\Phi$-Weyl bialternant $\chi_{_{\lambda}}$. 
In particular} $\WGF(L) = \chi_{_{\lambda}}${\sl , $L$ is rank symmetric and rank unimodal, and its rank-generating function} $\RGF(L,q)$ {\sl is given by equation (1) from the statement of \MainTheorem.}

Our goal here is to provide a completely uniform and from-scratch combinatorial proof of \MainTheorem. 
To do so, we apply a method of `vertex-coloring' that has been developed in some generality in  \cite{DonWeyl}. 
The idea is that in the presence of a suitably nice function that colors the vertices of a nicely structured poset, then a certain natural weight-generating function for the poset is automatically a specifiable Weyl bialternant or sum of Weyl bialternants. 
(For a precise statement of one version of this idea, see \ChainProductVersion\ below.) 
Using this methodology, our proof of the above main theorem only requires that we find a nice vertex-coloring function, see \GridPosetColoring.  

The vertex-coloring function presented in \GridPosetColoring\ was discovered by the authors empirically from a detailed analysis of the case-by-case descriptions of semistandard lattices as orderings of the semistandard tableaux presented in \cite{ADLMPPW}. 
But these tableaux descriptions can seem rather ad hoc and, when taken as a starting point, can be difficult to motivate for an uninitiated reader.  
We subsequently realized that our vertex-coloring function for semistandard lattices 
could be presented without reference to cases using only some very basic properties of two-color grid posets. 
This perspective serves to unify and greatly simplify some arguments and to distinguish those aspects of our approach that rely only on rudimentary, first-principles combinatorics. 

Our work in this paper leaves open the possibility that there might be other two-color grid posets and associated diamond-colored distributive lattices that meet the requirements of \MainTheorem. 
However, the second-listed author has shown that any two-color grid poset that meets the requirements of \MainTheorem\ must be one of the semistandard posets of \cite{ADLMPPW} (see Theorem 6.4 of \cite{DonDistributive}). 
That is, semistandard posets are the \underline{only} answers to the obvious combinatorial classification problem implicit in the statement of \MainTheorem. 

{\bf \S 2\ \ Some general background on posets and Weyl symmetric functions.} 
The combinatorial and Weyl-group theoretic conventions of this paper largely follow those of \cite{ADLP} and \cite{ADLMPPW}. 
In this section we highlight some of the key ideas necessary to understand the main results of the paper, and we state two general results about posets and Weyl bialternants / Weyl symmetric functions: \ChainProductVersion\ and \SplittingBasicsProposition. 
For a self-contained account of the foundations of Weyl symmetric function theory and Weyl bialternants, see \cite{DonWeyl}. 

Our starting point is a finite rank $n$ root system $\Phi$ residing in an $n$-dimensional Euclidean space $\EuclideanE$ with inner product $\langle \cdot,\cdot \rangle$. 
The related objects (coroots $\{\alpha^{\vee} = \frac{2}{\langle \alpha,\alpha \rangle}\alpha\}_{\alpha \in \Phi}$, simple roots $\{\alpha_{i}\}_{i \in I}$ over some index set $I$, Cartan matrix $M_{\Phi} = \left(\langle \alpha_{i},\alpha_{j}^{\vee} \rangle\right)_{i, j \in I}$, sets of positive and negative roots $\Phi^{+}$ and $\Phi^{-}$ respectively, fundamental weights $\{\omega_{i}\}_{i \in I}$ dual to the simple coroots, the lattice of weights $\Lambda = \{\sum_{i \in I}a_{i}\omega_{i}\, |\, a_{i} \in \mathbb{Z}, i \in I\}$, dominant weights $\Lambda^{+} = \{\sum_{i \in I}a_{i}\omega_{i}\, |\, a_{i} \in \mathbb{Z}_{\geq 0}, i \in I\}$, Weyl group $W$, longest Weyl group element $w_{0}$, special elements $\varrho := \sum \omega_{i} = \frac{1}{2}\sum_{\alpha \in \Phi^{+}}\alpha$ and $\varrho^{\vee} = \frac{1}{2}\sum_{\alpha \in \Phi^{+}}\alpha^{\vee}$, etc) are obtained as usual.  

The group ring $\mathbb{Z}[\Lambda]$ has as a $\mathbb{Z}$-basis the formal exponentials $\{e^{\mu}\}_{\lambda \in \Lambda}$.  
The Weyl group $W$ acts on $\mathbb{Z}[\Lambda]$ in the obvious way. 
The {\em ring of Weyl symmetric functions} $\mathbb{Z}[\Lambda]^{W}$ is the subring of $W$-invariants in $\mathbb{Z}[\Lambda]$.  
The {\em subgroup of $W$-alternants} $\mathbb{Z}[\Lambda]^{\mbox{\tiny alt}}$ consists of those group ring elements $\varphi$ for which $\sigma.\varphi = \det(\sigma)\varphi$ for all Weyl group elements $\sigma$. 
Define a mapping 
$\mathcal{A}: \mathbb{Z}[\lambda] \longrightarrow 
\mathbb{Z}[\Lambda]^{\mbox{\tiny alt}}$ 
by the rule $\mathcal{A}(\varphi) := \sum_{\sigma 
\in W}\det(\sigma) \sigma.\varphi$.  The Weyl denominator is the 
alternant $\mathcal{A}(e^{\varrho})$, which  factors as \[\mathcal{A}(e^{\varrho}) = e^{\varrho}\Big(\prod_{\alpha \in \Phi^{+}}(1-e^{-\alpha})\Big) = 
\prod_{\alpha \in \Phi^{+}}(e^{\alpha/2}-e^{-\alpha/2}) = 
e^{-\varrho}\Big(\prod_{\alpha 
\in \Phi^{+}}(e^{\alpha}-1)\Big).\]  
The following is a sort of fundamental theorem for Weyl symmetric functions: 
For any dominant weight $\lambda$, there exists a unique $\chi_{_{\lambda}} \in \mathbb{Z}[\Lambda]$ for which $\mathcal{A}(e^{\varrho})\, \chi_{_{\lambda}} = \mathcal{A}(e^{\lambda+\varrho})$, 
and in fact the $\chi_{_{\lambda}}$'s comprise a $\mathbb{Z}$-basis for $\mathbb{Z}[\Lambda]^{W}$. Each $\chi_{_{\lambda}}$ is a {\em Weyl bialternant}, sometimes also called a `Weyl character'.  
These are general Weyl group analogs of the Schur functions of classical symmetric function theory. 

The posets we work with here are finite. 
We identify a given poset $R$ with its Hasse diagram,  the directed graph whose edges depict the covering relations for the poset.  
In displayed figures of posets, edges are presented without arrowheads and are assumed to be directed upward. 
Also, many of the posets we work with are ranked. 
If $\rho: R \longrightarrow \{0, 1, 2, \ldots , l\}$ is a surjective rank function for a ranked poset $R$, then the number $l$ is the length of $R$ with respect to $\rho$, and the rank generating function for $R$ with respect to $\rho$ is the $q$-polynomial $\RGF(R,q) = \sum_{i = 0}^{l} |\rho^{-1}(i)|\, q^{i}$. 
A degree $l$ polynomial $a_{0}+a_{1}q+\cdots+a_{l}q^{l}$ in the variable $q$ is {\em symmetric} if $a_{i} = a_{l-i}$ for all $0 \leq i \leq l$ and is {\em unimodal} if for some $p$ ($0 \leq p \leq l$) it is the case that $a_{0} \leq \cdots \leq a_{p} \geq \cdots \geq a_{l}$. 
Then $R$ is {\em rank symmetric} (respectively, {\em rank unimodal}) if $\RGF(R,q)$ is a symmetric (resp.\ unimodal). 

For us, the posets that serve as models for Weyl bialternants are ranked and `edge-colored' by a function that assigns to each edge of $R$ a `color' from some set $I$.   
In this case we let $\rho_{i}(\xelt)$ denote the rank of an element $\xelt$ in $R$ within its $i$-component $\mathbf{comp}_{i}(\xelt)$ and $l_{i}(\xelt)$ denote the length of $\mathbf{comp}_{i}(\xelt)$. 
Set $\delta_i(\xelt) := l_i(\xelt)-\rho_i(\xelt)$. 

Now suppose $R$ is a ranked poset 
with edges colored by $I$, an index set of cardinality $n$ 
for a choice of simple roots for our  
finite root system $\Phi$.  
Let $wt: R \longrightarrow 
\Lambda$ be the function given by $wt(\xelt) = 
\sum_{i \in I}m_{i}(\xelt)\omega_{i}$ for all $\xelt \in R$, where 
$m_{i}(\xelt) := 2\rho_{i}(\xelt)-l_{i}(\xelt) = \rho_{i}(\xelt) - 
\delta_{i}(\xelt)$. Call $wt(\xelt)$ the 
{\em weight} of $\xelt$. 
Given the Cartan matrix $M_{\Phi}$ 
associated with our choice of simple 
roots for the root system $\Phi$, we say $R$ is 
$M_{\Phi}$-{\em structured} if it has 
the following property: 
$wt(\selt) + \alpha_{i} = 
wt(\telt)$ whenever $\selt \myarrow{i} \telt$ in $R$, or equivalently for 
all $j \in I$, $m_{j}(\selt) + \langle \alpha_{i},\alpha_{j}^{\vee} \rangle = m_{j}(\telt)$ 
whenever $\selt \myarrow{i} \telt$ in $R$. 

For such an $R$, the weight generating function $\WGF(R) := \displaystyle \sum_{\xelt \in R} e^{wt(\xelt)}$ is an 
element of the group ring $\mathbb{Z}[\Lambda]$. 
The $W$-invariance of $\WGF(R)$ is implied by a number 
of combinatorial conditions on $R$, such as the property that all of its $i$-components are rank symmetric, cf.\ Lemma 3.5 of \cite{DonWeyl}. 
We say $R$ is a {\em splitting poset for a Weyl symmetric function} $\chi$ if it is $M_{\Phi}$-structured and $\WGF(R) = \chi$. 

The next result is Corollary 8.2.B of \cite{DonWeyl}, whose proof does not depend upon the classification of finite root systems. 
To effect the statement of \ChainProductVersion, we need the following language and notation. 
Suppose $\mathcal{C}$ is a product $\mathcal{T}_{1} \times \cdots \times \mathcal{T}_{p}$ of chains $\mathcal{T}_{1}, \ldots ,\mathcal{T}_{p}$. 
A {\em face of} $\mathcal{C}$ is any subset of the form 
\[\left\{(\xelt_1,\ldots,\xelt_p) \in \mathcal{C}\, |\,  \mbox{$\xelt_q$ is maximal in $\mathcal{T}_{q}$}\right\},\]
where $q$ is some fixed index in $\{1, 2, \ldots, p\}$. 
A {\em sub-face of} $\mathcal{C}$ is a complement of a face. 
If $\phi: Q \longrightarrow \mathcal{C}$ is a poset isomorphism, then say $S \subseteq Q$ is a {\em sub-face of} $Q$ if $\phi(S)$ is a sub-face of $\mathcal{C}$. 
If the isomorphism $\phi$ is understood, we simply call $S$ a sub-face of $Q$. 

\noindent 
{\bf \ChainProductVersion}\ \ {\sl Suppose $R$ is an $M_{\Phi}$-structured poset. 
Suppose also that} $\WGF(R)$ {\sl is $W$-invariant. 
(The latter is guaranteed if, for example, for all $i \in I$ the $i$-components of $R$ are rank symmetric, cf.\ Lemma 3.5 of \cite{DonWeyl}.) 
Let $\mathcal{S} := \{\selt \in R\, |\, \delta_{i}(\selt) = 0 \mbox{ for all } i \in I\}$, so $wt(\selt) \in \Lambda^{+}$ for all $\selt \in \mathcal{S}$. 
Suppose $\kappa: R \setminus \mathcal{S} \longrightarrow I$ is a function such that for each $\xelt \in R \setminus \mathcal{S}$ with $i := \kappa(\xelt)$ we have (1) $\comp_{i}(\xelt)$ is isomorphic to a product of chains and (2) $\{\yelt \in \comp_{i}(\xelt)\, |\, \yelt \in R \setminus \mathcal{S} \mbox{ and } \kappa(\yelt) = i\}$ is a sub-face of $\comp_{i}(\xelt)$. Then $R$ is a splitting poset with} $\displaystyle \WGF(R) = \sum_{\selt \in \mathcal{S}}\chi_{_{wt(\selt)}}$.\hfill\QED

The next result follows from Proposition 4.7 of \cite{DonWeyl}. 
Unimodality of the rank-generating function is the only aspect of this result that depends upon the representation theory of semisimple Lie algebras, see for example Corollary 2.22 of \cite{DonWeyl}. 

\noindent 
{\bf \SplittingBasicsProposition}\ \ {\sl Let $\lambda \in \Lambda^{+}$ and suppose $R$ is a connected splitting poset for $\chi_{_{\lambda}}$. 
Then $R$ has a unique rank function. 
Moreover, $R$ is rank symmetric and rank unimodal, and its rank generating function is a polynomial of degree $2\langle \lambda,\varrho^{\vee} \rangle$ that can be written as follows, where the quantities $\langle \lambda+\varrho,\alpha^{\vee} \rangle$ and $\langle\varrho,\alpha^{\vee} \rangle$ are positive integers for each $\alpha \in \Phi^{+}$:} 

\hspace*{2in}$\displaystyle \RGF(R,q) = \mbox{\rule[-10mm]{0mm}{20mm}\large $\displaystyle \prod_{\alpha \in \Phi^{+}}$}\frac{1-q^{\langle \lambda+\varrho,\alpha^{\vee} \rangle}}{1-q^{\langle \varrho,\alpha^{\vee} \rangle}}$.\hfill\QED

Of course, one can specialize the preceding formula by letting $q \rightarrow 1$ in order to obtain a product-of-quotients expression for the cardinality of the splitting poset $R$. 

{\bf \S \TwoColorGridSection\ \ Two-color grid posets and their associated edge-colored distributive lattices.} 
Here we borrow liberally from the theory of two-color grid posets developed in \cite{ADLP} and \cite{ADLMPPW}.  
Let $m$ be a positive integer, and denote by $[m]$ the set $\{1, 2, \ldots , m\}$. 
Given a finite poset $(P,\leq_{_{P}})$, a {\em chain function for} $P$ is a function $\mathbf{chain}: P 
\longrightarrow [m]$ for some positive integer $m$ 
such that (1) $\mathbf{chain}^{-1}(i)$ is a (possibly empty) 
chain in $P$ for $1 \leq 
i \leq m$, and (2) given any cover  
$u \rightarrow v$ in $P$, it is the case that either 
$\mathbf{chain}(u) = \mathbf{chain}(v)$ or $\mathbf{chain}(u) = 
\mathbf{chain}(v) + 1$.  
A {\em grid poset} 
is a finite poset $(P,\leq_{_{P}})$ 
together with a chain function  
$\mathbf{chain}: P \longrightarrow [m]$ 
for some $m \geq 1$.  Depending on context, 
the notation $P$ can refer to 
the grid poset $(P,\leq_{_{P}}, 
\mathbf{chain}: P \longrightarrow [m])$ 
or the underlying poset $(P,\leq_{_{P}})$. 
The conditions on $\mathbf{chain}$ imply that an element in a grid 
poset covers no more than two elements and is covered by no more than 
two elements. 
{\sl Although it is not necessary for the chain function to be surjective, we will assume this is the case for the grid posets we consider here}, as it is convenient and costs us no generality. 
Then for $1 \leq i \leq m$ we set 
$\mathcal{C}_{i} := \mathbf{chain}^{-1}(i)$ and let $z_{i}$ denote the maximal element of $\mathcal{C}_{i}$. 
The {\em periphery} of $P$ is the set $\{z_{1}, z_{2}, \ldots , z_{m}\}$. 
When we depict grid 
posets, the chains $\mathcal{C}_{i}$ are directed from SW to NE.  
See \SomeFigures\ for examples. 

A {\em two-color function} for a grid poset $(P,\leq_{_{P}}, 
\mathbf{chain}: P \longrightarrow [m])$ is a function 
$\mathbf{color}: P \longrightarrow \Delta$ such that (1) $|\Delta| = 
2$, (2) $\mathbf{color}(u) = \color(v)$ 
if $\mathbf{chain}(u) = \mathbf{chain}(v)$, and (3) 
if $u$ and $v$ are in the same connected component 
of $P$ with $\mathbf{chain}(u) = \mathbf{chain}(v)+1$, 
then $\color(u) \not= \color(v)$. A {\em \dichromatic grid poset} is a 
grid poset 
$(P, \leq_{_{P}}, \mathbf{chain}: P \longrightarrow [m])$ 
together with a two-color function 
$\color: P \longrightarrow \Delta$. We sometimes use the 
notation $P$ by itself to refer to the \dichromatic grid poset 
$\digriddelta$. 
We typically use $\Delta = \{\alpha,\beta\}$. 

A two-color grid poset has the {\em max property} if (1) $\chain(u) \leq 2$ for any maximal element $u$ and (2) $\color(u) \not= \color(v)$ whenever $u$ and $v$ are distinct maximal elements.  
This innocuous property is one of the key hypotheses of \MainTheorem. 
It is one of two combinatorial properties that uniquely characterize semistandard posets, see \S 6 of \cite{DonDistributive}. 
Moreover, this property was crucial in formulating a uniform combinatorial argument for the uniqueness of weight bases for irreducible semisimple Lie algebra representations associated with semistandard lattices (when such weight bases exist), see \cite{ADLP}. 

For any two-color grid poset $P$, the associated distributive lattice $L = \mathbf{J}(P)$ of order ideals taken from $P$ is naturally edge-colored: When an order ideal $\telt$ covers an order ideal $\selt$ in the Hasse diagram for $L$, then $\telt \setminus \selt = \{u\}$ for some $u \in P$, in which case we give the edge $\selt \rightarrow \telt$ the color $i := \color(u)$ and depict this as $\selt \myarrow{i} \telt$. 
We write $L= \mathbf{J}_{color}(P)$ to emphasize that the edges of $L$ have colors thusly associated with the vertex colors of $P$.  
This distributive lattice is diamond-colored in the sense that on any diamond of edges 
\setlength{\unitlength}{0.4cm}
\begin{picture}(2.25,1.65)
\multiput(1,-0.5)(1,1){2}{\circle*{0.2}}
\multiput(0,0.5)(1,1){2}{\circle*{0.2}}
\multiput(1,-0.5)(1,1){2}{\line(-1,1){1}}
\multiput(1,-0.5)(-1,1){2}{\line(1,1){1}}
\put(0.3,0.85){\footnotesize $i$} \put(1.3,0.85){\footnotesize $j$}
\put(0.25,-0.15){\footnotesize $k$} \put(1.3,-0.15){\footnotesize $l$}
\end{picture}
in the order diagram for $L$, we have $i=l$ and $j=k$. 

Given such $P$ and $L$, fix $\telt \in L$ and $\gamma \in \Delta$.  
Let $\gamma'$ be the opposite color of $\gamma$, so $\{\gamma'\} = \Delta \setminus \{\gamma\}$. 
We claim that $\comp_{\gamma}(\telt)$, viewed as a subposet of $L$, is poset-isomorphic to a product of chains. 
Indeed, let $\mathcal{T}_{\gamma}^{(1)}, \ldots , \mathcal{T}_{\gamma}^{(k)}$ be the chains in $P$ of color $\gamma$, numbered so that $1 \leq i < j \leq k$ if and only if $\chain(w^{(i)}) < \chain(w^{(j)})$ when $w^{(i)}$ and $w^{(j)}$ are the periphery elements of $\mathcal{T}_{\gamma}^{(i)}$ and $\mathcal{T}_{\gamma}^{(j)}$. 
Observe that for any $\selt \in \comp_{\gamma}(\telt)$, the color $\gamma'$ vertices of $\selt$ and $\telt$ are exactly the same. 
Now each pair of color $\gamma$ chains in $P$ must have a color $\gamma'$ chain between them, or else the grid poset disconnects between these two chains. 
Either way, it follows that as we maneuver through $\comp_{\gamma}(\telt)$, vertices from a given color $\gamma$ chain $\mathcal{T}_{\gamma}^{(i)}$ can be added to or deleted from $\telt$ independently of making any additions/deletions of vertices from a different chain $\mathcal{T}_{\gamma}^{(j)}$. 
In this way, we can discern that $\comp_{\gamma}(\telt)$ is isomorphic to a product of chains. 

We can formulate this precisely as follows. 
Pick $j \in \{1, 2, \ldots , k\}$. 
Suppose that for subsets $\mathcal{J}$ and $\mathcal{J}'$ of $\mathcal{T}_{\gamma}^{(j)}$ it is that case that each of $\telt \cup \mathcal{J}$, $\telt \setminus \mathcal{J}$, $\telt \cup \mathcal{J}'$, and $\telt \setminus \mathcal{J}'$ is an order ideal. 
Set $\mathcal{I} := \mathcal{J} \cup \mathcal{J}'$. 
It is easy to see that $\telt \cup \mathcal{I}$ and $\telt \setminus \mathcal{I}$ are also order ideals. 
So there is a largest subset of $\mathcal{T}_{\gamma}^{(j)}$, which we denote $\mathcal{I}_{j}$, for which $\telt \cup \mathcal{I}_{j}$ and $\telt \setminus \mathcal{I}_{j}$ are order ideals. 
In fact, $\mathcal{I}_{j}$ is an interval (possibly empty) within the chain $\mathcal{T}_{\gamma}^{(j)}$:  Suppose $x \leq v \leq y$ in $\mathcal{T}_{\gamma}^{(j)}$ with $x$ and $y$ in $\mathcal{I}_{j}$. 
Since $\telt \cup \mathcal{I}_{j}$ is an order ideal containing $y$ and since $v \leq y$, it follows that $v \in \telt \cup \mathcal{I}_{j}$. 
If $v \not\in \mathcal{I}_{j}$, then $v$ is in the order ideal $\telt \setminus \mathcal{I}_{j}$, so the fact that $x \leq v$ means that $x$ is in $\telt \setminus \mathcal{I}_{j}$ as well. 
But this contradicts the fact that $x$ was chosen from $\mathcal{I}_{j}$. 
We conclude that we must have $v \in \mathcal{I}_{j}$ as well. 
So $\mathcal{I}_{j}$ is an interval in the chain $\mathcal{T}_{\gamma}^{(j)}$, as claimed. 
It follows that $\mathbf{J}(\mathcal{I}_{j})$ is also a chain. 
Let $\phi: \comp_{\gamma}(\telt) \longrightarrow \mathbf{J}(\mathcal{I}_{1}) \times \cdots \times \mathbf{J}(\mathcal{I}_{k})$ be given by $\phi(\selt) = (\selt \cap \mathcal{I}_{1}, \ldots , \selt \cap \mathcal{I}_{k})$. 
It is routine to check that $\phi$ is a poset isomorphism, a fact we record as follows. 

\noindent 
{\bf \ComponentsAreChains}\ \ {\sl In the notation of the preceding paragraph, the function $\phi: \comp_{\gamma}(\telt) \longrightarrow \mathbf{J}(\mathcal{I}_{1}) \times \cdots \times \mathbf{J}(\mathcal{I}_{k})$ is a well-defined poset isomorphism.  
In particular, $\comp_{\gamma}(\telt)$ is isomorphic to a product of chains.}\hfill\QED

For more on the following discussion of `decomposing' grid posets and 
two-color grid posets, see \cite{ADLMPPW}.  
Let  $P$ be a grid poset with chain function $\mathbf{chain}: P 
\longrightarrow \{1,2,\ldots,m\}$.  Suppose $P_{1}$ is 
a nonempty order ideal and is a proper subset of $P$.  
Regard $P_{1}$ and $P_{2} := P \setminus P_{1}$ to be subposets of 
the poset $P$ in the induced order.  Suppose that whenever $u$ is a maximal 
(respectively minimal) element of $P_{1}$ and $v$ is a maximal 
(respectively minimal) element 
of $P_{2}$, then $\mathbf{chain}(u) \leq \mathbf{chain}(v)$.  Then we 
say that $P$ 
{\em decomposes into} $P_{1} \triangleleft P_{2}$, and we write 
$P = P_{1} \triangleleft P_{2}$. 
If no such order ideal $P_{1}$ exists, 
then we say the grid poset $P$ is {\em indecomposable}. 
If $P$ is a grid poset that decomposes into $P_{1} \triangleleft Q$, 
and if $Q$ decomposes into $P_{2} \triangleleft P_{3}$, then  
$P = P_{1} \triangleleft (P_{2} \triangleleft P_{3})$.  But now 
observe that $P = (P_{1} \triangleleft P_{2}) \triangleleft P_{3}$.  
So we may write $P = P_{1} \triangleleft P_{2} \triangleleft P_{3}$ 
unambiguously.  In general, if 
$P = P_{1} \triangleleft P_{2} \triangleleft \cdots \triangleleft 
P_{k}$, then each $P_{i}$ with chain function 
$\mathbf{chain}|_{P_{i}}$ is a grid subposet of $P$. 
If in addition $P$ is a \dichromatic grid poset with two-color 
function $\color$, then each $P_{i}$ with chain function 
$\mathbf{chain}|_{P_{i}}$ and two-color function 
$\color|_{P_{i}}$ is a \dichromatic grid subposet of $P$, and so 
$P_{1} \triangleleft P_{2} \triangleleft 
\cdots \triangleleft P_{k}$ is a 
decomposition of $P$ into \dichromatic grid posets.

{\bf \S \SemistandardSection\ \ Semistandard posets and lattices.} 
Let $\Phi$ denote 
a rank two root system with simple root basis $\{\alpha,\beta\}$, where $\alpha$ is short.  
The vertex colors and edge colors for the posets and lattices 
we now present are simple roots.  
Let $\omega_{\alpha} = \omega_{1} = (1,0)$ and 
$\omega_{\beta} = \omega_{2} = (0,1)$ respectively denote the 
corresponding fundamental weights. 
Any weight $\mu$ in $\Lambda$ of the form $\mu = p\omega_{\alpha} + 
q\omega_{\beta}$ (where $p$ and $q$ are integers) is now identified 
with the pair $(p,q)$ in $\mathbb{Z} \times \mathbb{Z}$. 
Then $\alpha$ and $\beta$ are respectively identified with the first 
and second 
row vectors from the Cartan matrix $M_\Phi$, cf.\ \S \StartSection.  

With respect to this notation, we define the $\Phi$-{\em 
fundamental posets} $P_{\Phi}(1,0)$ and 
$P_{\Phi}(0,1)$ to be the  
\dichromatic grid posets of \FundPosets. 
The corresponding $\Phi$-{\em fundamental 
lattices} are the edge-colored 
lattices $L_{\Phi}(1,0) := 
\mathbf{J}_{color}(P_{\Phi}(1,0))$ 
and 
$L_{\Phi}(0,1) := 
\mathbf{J}_{color}(P_{\Phi}(0,1))$ 
respectively.  Now let 
$\lambda = (a,b)$ be a pair of nonnegative integers.  There 
are exactly two possible ways that a \dichromatic grid poset $P$ 
with the max property 
can decompose as $P_{1} \triangleleft P_{2} \triangleleft \cdots \triangleleft 
P_{a+b}$ with $a$ of the $P_{i}$'s vertex-color isomorphic to  
$P_{\Phi}(1,0)$ and the remaining $P_{i}$'s vertex-color 
isomorphic to $P_{\Phi}(0,1)$: we will either have $P_{i}$ 
isomorphic to $P_{\Phi}(0,1)$ for $1 \leq i \leq b$ and 
isomorphic to $P_{\Phi}(1,0)$ for $1+b \leq i \leq a+b$ 
(in which 
case we set $P_{\Phi}^{\beta\alpha}(\lambda) := P$), 
or we will have $P_{i}$ 
isomorphic to $P_{\Phi}(1,0)$ for $1 \leq i \leq a$ and 
isomorphic to $P_{\Phi}(0,1)$ for $a+1 \leq i \leq a+b$ 
(in which 
case we set $P_{\Phi}^{\alpha\beta}(\lambda) := P$).  
Note that 
$P_{\Phi}^{\beta\alpha}(1,0) = 
P_{\Phi}^{\alpha\beta}(1,0) = 
P_{\Phi}(1,0)$, and   
$P_{\Phi}^{\beta\alpha}(0,1) = 
P_{\Phi}^{\alpha\beta}(0,1) = 
P_{\Phi}(0,1)$. When $a = b = 0$, then 
$P_{\Phi}^{\beta\alpha}(\lambda)$ and 
$P_{\Phi}^{\alpha\beta}(\lambda)$ are the empty set. 
We call $P_{\Phi}^{\beta\alpha}(\lambda)$ and 
$P_{\Phi}^{\alpha\beta}(\lambda)$ the 
$\Phi$-{\em semistandard posets} associated 
to $\lambda$.  For each rank two root system $\Phi$, 
$P_{\Phi}^{\beta\alpha}(2,2)$ is depicted in \GridPosets; 
$P_{\Phi}^{\alpha\beta}(2,2)$ is depicted in \GridPosetsII. 
The $\Phi$-{\em semistandard lattices} associated to 
$\lambda$ are the edge-colored lattices 
$\Lba := 
\mathbf{J}_{color}(P_{\Phi}^{\beta\alpha}(\lambda))$ 
and $L_{\Phi}^{\alpha\beta}(\lambda) := 
\mathbf{J}_{color}(P_{\Phi}^{\alpha\beta}(\lambda))$. 
Note that 
$L_{\Phi}^{\beta\alpha}(1,0) = 
L_{\Phi}^{\alpha\beta}(1,0) = 
L_{\Phi}(1,0)$, and   
$L_{\Phi}^{\beta\alpha}(0,1) = 
L_{\Phi}^{\alpha\beta}(0,1) = 
L_{\Phi}(0,1)$.

Proposition 4.2 of \cite{ADLMPPW} shows, using mostly uniform arguments concerning two-color grid posets, that any $\Phi$-semistandard lattice is $M_{\Phi}$-structured. 
It follows from the definitions that the above $\Phi$-semistandard posets have the max property, as claimed in the proposition below. 

\noindent
{\bf \StructureMaxResult}\ \ {\sl Let 
$\lambda = (a,b)$ be a pair of nonnegative integers, and let $L$ be 
one of the 
$\Phi$-semistandard lattices $L_{\Phi}^{\beta\alpha}(\lambda)$ 
or $L_{\Phi}^{\alpha\beta}(\lambda)$.  
Let $\selt \myarrow{\gamma} \telt$ be an 
edge of color $\gamma \in \{\alpha,\beta\}$ in $L$.  Then 
$wt(\selt) + \gamma = wt(\telt)$, and hence $L$ is $M_\Phi$-structured. 
The associated $\Phi$-semistandard posets $P_{\Phi}^{\beta\alpha}(\lambda)$ and 
$P_{\Phi}^{\alpha\beta}(\lambda)$ have the max property.}\hfill\QED

\newpage 
\begin{center}
\FundPosets: Fundamental posets for rank two root systems

\

\vspace*{-0.15in} 

\begin{tabular}{|c||c|c|}
\hline
Root system $\Phi$ & $P_{\Phi}(1,0)$ & 
$P_{\Phi}(0,1)$\\
\hline
\hline
\parbox[b]{0.5in}{\ 

$\myA_{1} \oplus \myA_{1}$

\ 

\vspace*{-0.1in}
} & \AOneAOneAlpha & \AOneAOneBeta\\
\hline
\parbox[b]{0.225in}{$\myA_{2}$
\vspace*{0.4in}} 
 & \ATwoAlpha & \ATwoBeta\\
\hline
\parbox[b]{0.225in}{$\myC_{2}$
\vspace*{1.1in}} 
 & \BTwoAAlpha & \BTwoBBeta\\
\hline
\parbox[b]{0.225in}{$\myG_{2}$
\vspace*{2.5in}} 
 & \GTwoAAlpha & \GTwoBBeta\\
\hline
\end{tabular}
\end{center}

\newpage
\begin{center}
\setlength{\unitlength}{1.1cm}
\begin{picture}(14,18.75)
\put(0.105,18.4){\FundLatticeIdealsFigure: Elements of fundamental  
lattices as order ideals taken from fundamental posets.}
\put(0.025,17.6){\parbox{6.15in}{\small \begin{center}
(Each order ideal taken from the fundamental poset is identified by the 
indices of its maximal vertices. 
For example, $\langle 6,7 \rangle$ in $L_{\mysmallG_{2}}(0,1)$ 
denotes the order ideal $\{v_{6}, v_{7}, v_{8}, v_{9}, v_{10}\}$ 
taken from $P_{\mysmallG_{2}}(0,1)$.) 
\end{center}
}}
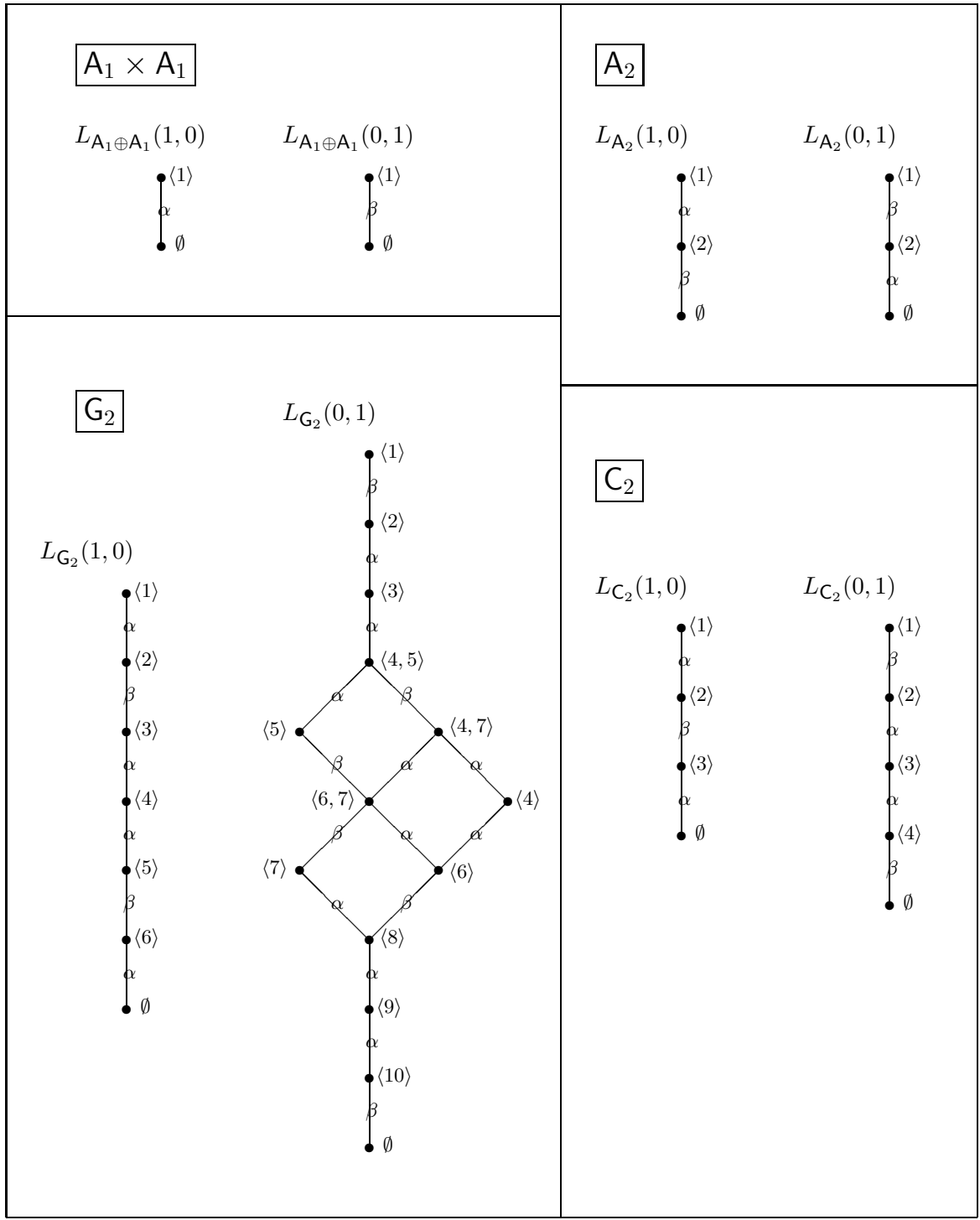
\put(0,-0.5){
\begin{picture}(14,18.5)
\put(0,0){\line(0,1){17.5}}
\put(8,0){\line(0,1){17.5}}
\put(14,0){\line(0,1){17.5}}
\put(0,0){\line(1,0){14}}
\put(0,13){\line(1,0){8}}
\put(8,12){\line(1,0){6}}
\put(0,17.5){\line(1,0){14}}
\put(1,16.5){\fbox{\Large $\myA_{1} \times \myA_{1}$}}
\put(0.5,13){\AoneAlphaIdeals}
\put(3.5,13){\AoneBetaIdeals}
\put(8.5,16.5){\fbox{\Large $\myA_{2}$}}
\put(8,12){\AtwoAlphaIdeals}
\put(11,12){\AtwoBetaIdeals}
\put(8.5,10.5){\fbox{\Large $\myC_{2}$}}
\put(8,4.5){\BtwoAlphaIdeals}
\put(11,3.5){\BtwoBetaIdeals}
\put(1,11.5){\fbox{\Large $\myG_{2}$}}
\put(0,2){\GtwoAlphaIdeals}
\put(2,0){\GtwoBetaIdeals}
\end{picture}
}
\end{picture}
\end{center} 


\newpage
\begin{center}
\GridPosets: Depicted below are four \dichromatic grid posets each 
possessing the max 
property.\\  
{\small (Each is a $\Phi$-semistandard poset 
$P_{\Phi}^{\beta\alpha}(2,2)$ as in \S \SemistandardSection; periphery vertices are indicated.)} 

\setlength{\unitlength}{1cm}
\begin{picture}(6,7)
\put(0,7){\fbox{$\Phi = \myA_{1} \oplus \myA_{1}$}}
\put(0.5,2.75){
\begin{picture}(3,3.5)
\put(3,2){\circle*{0.15}} 
\put(3.2,1.9){\footnotesize $\alpha$}
\put(4,3){\circle*{0.15}} 
\put(3.4,2.9){\footnotesize $z_{2}$}
\put(4.2,2.9){\footnotesize $\alpha$}
\put(0,2){\circle*{0.15}} 
\put(0.2,1.9){\footnotesize $\beta$}
\put(1,3){\circle*{0.15}} 
\put(0.4,2.9){\footnotesize $z_{1}$}
\put(1.2,2.9){\footnotesize $\beta$}
\put(3,2){\line(1,1){1}}
\put(0,2){\line(1,1){1}} 
\put(1.5,3.5){$\mathcal{C}_{1}$}
\put(4.5,3.5){$\mathcal{C}_{2}$}
\put(1.8,2.375){$\bigoplus$}
\end{picture}
}
\put(0,2){\fbox{$\Phi = \myA_{2}$}}
\put(0,-2.5){
\begin{picture}(3,3.5)
\put(4,2){\circle*{0.15}} 
\put(4.2,1.9){\footnotesize $\beta$}
\put(5,3){\circle*{0.15}} 
\put(4.4,2.9){\footnotesize $z_{3}$}
\put(5.2,2.9){\footnotesize $\beta$}
\put(1,1){\circle*{0.15}} 
\put(1.2,0.9){\footnotesize $\alpha$}
\put(2,2){\circle*{0.15}} 
\put(2.2,1.9){\footnotesize $\alpha$}
\put(3,3){\circle*{0.15}} 
\put(3.2,2.9){\footnotesize $\alpha$}
\put(4,4){\circle*{0.15}} 
\put(3.4,3.9){\footnotesize $z_{2}$}
\put(4.2,3.9){\footnotesize $\alpha$}
\put(0,2){\circle*{0.15}} 
\put(0.2,1.9){\footnotesize $\beta$}
\put(1,3){\circle*{0.15}} 
\put(0.4,2.9){\footnotesize $z_{1}$}
\put(1.2,2.9){\footnotesize $\beta$}
\put(1,1){\line(1,1){3}}
\put(4,2){\line(1,1){1}}
\put(0,2){\line(1,1){1}} 
\put(1,1){\line(-1,1){1}} 
\put(2,2){\line(-1,1){1}} 
\put(3,3){\line(1,-1){1}}
\put(4,4){\line(1,-1){1}}
\put(1.5,3.5){$\mathcal{C}_{1}$}
\put(4.5,4.5){$\mathcal{C}_{2}$}
\put(5.5,3.5){$\mathcal{C}_{3}$}
\end{picture}
}
\end{picture}
\hspace*{1.5cm}
\setlength{\unitlength}{1cm}
\begin{picture}(8,8)
\put(0,6){\fbox{$\Phi = \myC_{2}$}}
\put(0,0.25){
\begin{picture}(3,3.5)
\put(6,3){\circle*{0.15}} 
\put(6.2,2.9){\footnotesize $\alpha$}
\put(7,4){\circle*{0.15}} 
\put(6.4,3.9){\footnotesize $z_{4}$}
\put(7.2,3.9){\footnotesize $\alpha$}
\put(1,0){\circle*{0.15}} 
\put(1.2,-0.1){\footnotesize $\beta$}
\put(3,2){\circle*{0.15}} 
\put(3.2,1.9){\footnotesize $\beta$}
\put(5,4){\circle*{0.15}} 
\put(5.2,3.9){\footnotesize $\beta$}
\put(6,5){\circle*{0.15}} 
\put(5.4,4.9){\footnotesize $z_{3}$}
\put(6.2,4.9){\footnotesize $\beta$}
\put(0,1){\circle*{0.15}} 
\put(0.2,0.9){\footnotesize $\alpha$}
\put(1,2){\circle*{0.15}} 
\put(1.2,1.9){\footnotesize $\alpha$}
\put(2,3){\circle*{0.15}} 
\put(2.2,2.9){\footnotesize $\alpha$}
\put(3,4){\circle*{0.15}} 
\put(3.2,3.9){\footnotesize $\alpha$}
\put(4,5){\circle*{0.15}} 
\put(4.2,4.9){\footnotesize $\alpha$}
\put(5,6){\circle*{0.15}} 
\put(4.4,5.9){\footnotesize $z_{2}$}
\put(5.2,5.9){\footnotesize $\alpha$}
\put(0,3){\circle*{0.15}} 
\put(0.2,2.9){\footnotesize $\beta$}
\put(2,5){\circle*{0.15}} 
\put(1.4,4.9){\footnotesize $z_{1}$}
\put(2.2,4.9){\footnotesize $\beta$}
\put(0,1){\line(1,1){5}}
\put(1,0){\line(1,1){5}}
\put(0,3){\line(1,1){2}} 
\put(6,3){\line(1,1){1}} 
\put(1,0){\line(-1,1){1}} 
\put(1,2){\line(-1,1){1}} 
\put(3,2){\line(-1,1){1}} 
\put(3,4){\line(-1,1){1}} 
\put(4,5){\line(1,-1){2}}
\put(5,6){\line(1,-1){2}}
\put(2.5,5.5){$\mathcal{C}_{1}$}
\put(5.5,6.5){$\mathcal{C}_{2}$}
\put(6.5,5.5){$\mathcal{C}_{3}$}
\put(7.5,4.5){$\mathcal{C}_{4}$}
\end{picture}
}
\end{picture}

\setlength{\unitlength}{1cm}
\begin{picture}(11,13)
\put(-0.5,9){\fbox{$\Phi = \myG_{2}$}}
\put(2,0){\circle*{0.15}} 
\put(2.2,-0.15){\footnotesize $\beta$}
\put(1,1){\circle*{0.15}} 
\put(1.2,0.9){\footnotesize $\alpha$}
\put(2,2){\circle*{0.15}} 
\put(2.2,1.9){\footnotesize $\alpha$}
\put(1,3){\circle*{0.15}} 
\put(1.2,2.9){\footnotesize $\beta$}
\put(3,3){\circle*{0.15}} 
\put(3.2,2.9){\footnotesize $\alpha$}
\put(5,3){\circle*{0.15}} 
\put(5.2,2.85){\footnotesize $\beta$}
\put(0,4){\circle*{0.15}} 
\put(0.2,3.9){\footnotesize $\alpha$}
\put(2,4){\circle*{0.15}} 
\put(2.2,3.9){\footnotesize $\beta$}
\put(4,4){\circle*{0.15}} 
\put(4.2,3.9){\footnotesize $\alpha$}
\put(1,5){\circle*{0.15}} 
\put(1.2,4.9){\footnotesize $\alpha$}
\put(5,5){\circle*{0.15}} 
\put(5.2,4.9){\footnotesize $\alpha$}
\put(9,5){\circle*{0.15}} 
\put(9.2,4.9){\footnotesize $\alpha$}
\put(2,6){\circle*{0.15}} 
\put(2.2,5.9){\footnotesize $\alpha$}
\put(4,6){\circle*{0.15}} 
\put(4.2,5.9){\footnotesize $\beta$}
\put(6,6){\circle*{0.15}} 
\put(6.2,5.9){\footnotesize $\alpha$}
\put(8,6){\circle*{0.15}} 
\put(8.2,5.9){\footnotesize $\beta$}
\put(1,7){\circle*{0.15}} 
\put(1.2,6.9){\footnotesize $\beta$}
\put(3,7){\circle*{0.15}} 
\put(3.2,6.9){\footnotesize $\alpha$}
\put(5,7){\circle*{0.15}} 
\put(5.2,6.9){\footnotesize $\beta$}
\put(7,7){\circle*{0.15}} 
\put(7.2,6.9){\footnotesize $\alpha$}
\put(11,7){\circle*{0.15}} 
\put(10.4,6.9){\footnotesize $z_{6}$}
\put(11.2,6.9){\footnotesize $\alpha$}
\put(4,8){\circle*{0.15}} 
\put(4.2,7.9){\footnotesize $\alpha$}
\put(8,8){\circle*{0.15}} 
\put(8.2,7.9){\footnotesize $\alpha$}
\put(10,8){\circle*{0.15}} 
\put(9.4,7.9){\footnotesize $z_{5}$}
\put(10.2,7.9){\footnotesize $\beta$}
\put(5,9){\circle*{0.15}} 
\put(5.2,8.9){\footnotesize $\alpha$}
\put(7,9){\circle*{0.15}} 
\put(7.2,8.9){\footnotesize $\beta$}
\put(9,9){\circle*{0.15}} 
\put(9.2,8.9){\footnotesize $\alpha$}
\put(4,10){\circle*{0.15}} 
\put(3.4,9.9){\footnotesize $z_{1}$}
\put(4.2,9.9){\footnotesize $\beta$}
\put(6,10){\circle*{0.15}} 
\put(6.2,9.9){\footnotesize $\alpha$}
\put(10,10){\circle*{0.15}} 
\put(9.4,9.9){\footnotesize $z_{4}$}
\put(10.2,9.9){\footnotesize $\alpha$}
\put(9,11){\circle*{0.15}} 
\put(8.4,10.9){\footnotesize $z_{3}$}
\put(9.2,10.9){\footnotesize $\beta$}
\put(8,12){\circle*{0.15}} 
\put(7.4,11.9){\footnotesize $z_{2}$}
\put(8.2,11.9){\footnotesize $\alpha$}
\put(9,5){\line(1,1){2}} 
\put(2,0){\line(1,1){8}} 
\put(1,1){\line(1,1){9}} 
\put(1,3){\line(1,1){8}} 
\put(0,4){\line(1,1){8}} 
\put(1,7){\line(1,1){3}}
\put(2,0){\line(-1,1){1}} 
\put(2,2){\line(-1,1){2}} 
\put(3,3){\line(-1,1){2}} 
\put(5,3){\line(-1,1){1}} 
\put(2,6){\line(-1,1){1}} 
\put(5,5){\line(-1,1){2}} 
\put(6,6){\line(-1,1){2}} 
\put(9,5){\line(-1,1){2}} 
\put(5,9){\line(-1,1){1}} 
\put(8,8){\line(-1,1){2}} 
\put(11,7){\line(-1,1){2}} 
\put(10,10){\line(-1,1){2}} 
\put(11.5,7.5){$\mathcal{C}_{6}$}
\put(10.5,8.5){$\mathcal{C}_{5}$}
\put(10.5,10.5){$\mathcal{C}_{4}$}
\put(9.5,11.5){$\mathcal{C}_{3}$}
\put(8.5,12.5){$\mathcal{C}_{2}$}
\put(4.5,10.5){$\mathcal{C}_{1}$}
\end{picture} 
\end{center}

\newpage
\begin{center}
\GridPosetsII: Depicted below are four \dichromatic grid posets each 
possessing the max 
property.\\  
{\small (Each is a $\Phi$-semistandard poset 
$P_{\Phi}^{\alpha\beta}(2,2)$ as in \S \SemistandardSection; periphery vertices are indicated.)} 

\setlength{\unitlength}{1cm}
\begin{picture}(6,7)
\put(0,7){\fbox{$\Phi = \myA_{1} \oplus \myA_{1}$}}
\put(0.5,2.75){
\begin{picture}(3,3.5)
\put(3,2){\circle*{0.15}} 
\put(3.2,1.9){\footnotesize $\beta$}
\put(4,3){\circle*{0.15}} 
\put(3.4,2.9){\footnotesize $z_{2}$}
\put(4.2,2.9){\footnotesize $\beta$}
\put(0,2){\circle*{0.15}} 
\put(0.2,1.9){\footnotesize $\alpha$}
\put(1,3){\circle*{0.15}} 
\put(0.4,2.9){\footnotesize $z_{1}$}
\put(1.2,2.9){\footnotesize $\alpha$}
\put(3,2){\line(1,1){1}}
\put(0,2){\line(1,1){1}} 
\put(1.5,3.5){$\mathcal{C}_{1}$}
\put(4.5,3.5){$\mathcal{C}_{2}$}
\put(1.8,2.375){$\bigoplus$}
\end{picture}
}
\put(0,2){\fbox{$\Phi = \myA_{2}$}}
\put(0,-2.5){
\begin{picture}(3,3.5)
\put(4,2){\circle*{0.15}} 
\put(4.2,1.9){\footnotesize $\alpha$}
\put(5,3){\circle*{0.15}} 
\put(4.4,2.9){\footnotesize $z_{3}$}
\put(5.2,2.9){\footnotesize $\alpha$}
\put(1,1){\circle*{0.15}} 
\put(1.2,0.9){\footnotesize $\beta$}
\put(2,2){\circle*{0.15}} 
\put(2.2,1.9){\footnotesize $\beta$}
\put(3,3){\circle*{0.15}} 
\put(3.2,2.9){\footnotesize $\beta$}
\put(4,4){\circle*{0.15}} 
\put(3.4,3.9){\footnotesize $z_{2}$}
\put(4.2,3.9){\footnotesize $\beta$}
\put(0,2){\circle*{0.15}} 
\put(0.2,1.9){\footnotesize $\alpha$}
\put(1,3){\circle*{0.15}} 
\put(0.4,2.9){\footnotesize $z_{1}$}
\put(1.2,2.9){\footnotesize $\alpha$}
\put(1,1){\line(1,1){3}}
\put(4,2){\line(1,1){1}}
\put(0,2){\line(1,1){1}} 
\put(1,1){\line(-1,1){1}} 
\put(2,2){\line(-1,1){1}} 
\put(3,3){\line(1,-1){1}}
\put(4,4){\line(1,-1){1}}
\put(1.5,3.5){$\mathcal{C}_{1}$}
\put(4.5,4.5){$\mathcal{C}_{2}$}
\put(5.5,3.5){$\mathcal{C}_{3}$}
\end{picture}
}
\end{picture}
\hspace*{1.5cm}
\setlength{\unitlength}{1cm}
\begin{picture}(8,8)
\put(0,6){\fbox{$\Phi = \myC_{2}$}}
\put(0,0.25){
\begin{picture}(3,3.5)
\put(1,3){\ZVertexForPosets{1}{\alpha}}
\put(0,2){\VertexForPosets{2}{\alpha}}
\put(6,6){\ZVertexForPosets{2}{\beta}}
\put(4,4){\VertexForPosets{4}{\beta}}
\put(2,2){\VertexForPosets{5}{\beta}}
\put(1,1){\VertexForPosets{6}{\beta}}
\put(7,5){\ZVertexForPosets{3}{\alpha}}
\put(6,4){\VertexForPosets{8}{\alpha}}
\put(5,3){\VertexForPosets{9}{\alpha}}
\put(4,2){\VertexForPosets{10}{\alpha}}
\put(3,1){\VertexForPosets{11}{\alpha}}
\put(2,0){\VertexForPosets{12}{\alpha}}
\put(7,3){\ZVertexForPosets{4}{\beta}}
\put(5,1){\VertexForPosets{14}{\beta}}
\put(0,2){\line(1,1){1}}
\put(1,1){\line(1,1){5}}
\put(2,0){\line(1,1){5}} 
\put(5,1){\line(1,1){2}} 
\put(2,0){\line(-1,1){2}} 
\put(3,1){\line(-1,1){2}} 
\put(5,1){\line(-1,1){1}} 
\put(5,3){\line(-1,1){1}} 
\put(7,3){\line(-1,1){1}}
\put(7,5){\line(-1,1){1}}
\put(1.5,3.5){$\mathcal{C}_{1}$}
\put(6.5,6.5){$\mathcal{C}_{2}$}
\put(7.5,5.5){$\mathcal{C}_{3}$}
\put(7.5,3.5){$\mathcal{C}_{4}$}
\end{picture}
}
\end{picture}

\setlength{\unitlength}{1cm}
\begin{picture}(11,13)
\put(-0.5,9){\fbox{$\Phi = \myG_{2}$}}
\put(2,7){\ZVertexForPosets{1}{\alpha}}
\put(0,5){\VertexForPosets{2}{\alpha}}
\put(9,12){\ZVertexForPosets{2}{\beta}}
\put(6,9){\VertexForPosets{4}{\beta}}
\put(3,6){\VertexForPosets{5}{\beta}}
\put(1,4){\VertexForPosets{6}{\beta}}
\put(10,11){\ZVertexForPosets{3}{\alpha}}
\put(9,10){\VertexForPosets{8}{\alpha}}
\put(8,9){\VertexForPosets{9}{\alpha}}
\put(7,8){\VertexForPosets{10}{\alpha}}
\put(6,7){\VertexForPosets{11}{\alpha}}
\put(5,6){\VertexForPosets{12}{\alpha}}
\put(4,5){\VertexForPosets{13}{\alpha}}
\put(3,4){\VertexForPosets{14}{\alpha}}
\put(2,3){\VertexForPosets{15}{\alpha}}
\put(1,2){\VertexForPosets{16}{\alpha}}
\put(10,9){\ZVertexForPosets{4}{\beta}}
\put(9,8){\VertexForPosets{18}{\beta}}
\put(7,6){\VertexForPosets{19}{\beta}}
\put(6,5){\VertexForPosets{20}{\beta}}
\put(4,3){\VertexForPosets{21}{\beta}}
\put(2,1){\VertexForPosets{22}{\beta}}
\put(11,8){\ZVertexForPosets{5}{\alpha}}
\put(10,7){\VertexForPosets{24}{\alpha}}
\put(9,6){\VertexForPosets{25}{\alpha}}
\put(8,5){\VertexForPosets{26}{\alpha}}
\put(7,4){\VertexForPosets{27}{\alpha}}
\put(6,3){\VertexForPosets{28}{\alpha}}
\put(5,2){\VertexForPosets{29}{\alpha}}
\put(3,0){\VertexForPosets{30}{\alpha}}
\put(10,5){\ZVertexForPosets{6}{\beta}}
\put(7,2){\VertexForPosets{32}{\beta}}
\put(7,2){\line(1,1){3}} 
\put(3,0){\line(1,1){8}} 
\put(2,1){\line(1,1){8}} 
\put(1,2){\line(1,1){9}} 
\put(1,4){\line(1,1){8}} 
\put(0,5){\line(1,1){2}}
\put(3,0){\line(-1,1){2}} 
\put(2,3){\line(-1,1){2}} 
\put(5,2){\line(-1,1){2}} 
\put(4,5){\line(-1,1){2}} 
\put(7,2){\line(-1,1){1}} 
\put(7,4){\line(-1,1){2}} 
\put(8,5){\line(-1,1){2}} 
\put(7,8){\line(-1,1){1}} 
\put(10,5){\line(-1,1){1}} 
\put(10,7){\line(-1,1){2}} 
\put(11,8){\line(-1,1){2}} 
\put(10,11){\line(-1,1){1}} 
\put(10.5,5.5){$\mathcal{C}_{6}$}
\put(11.5,8.5){$\mathcal{C}_{5}$}
\put(10.5,9.5){$\mathcal{C}_{4}$}
\put(10.5,11.5){$\mathcal{C}_{3}$}
\put(9.5,12.5){$\mathcal{C}_{2}$}
\put(2.5,7.5){$\mathcal{C}_{1}$}
\end{picture} 
\end{center}


\newpage 
{\bf \S \ColorFunctionSection\ \ Vertex-coloring.} 
The main result of this section requires some further combinatorial set-up. 
In particular, we must define a function that `colors' the vertices of the diamond-colored distributive lattice associated with a two-color grid poset. 
For the remainder of the section, $P$ is a two-color grid poset with surjective chain function $\chain: P \longrightarrow [m]$, two-color function $\color: P \longrightarrow \Delta$, and associated edge-colored distributive lattice $L = \mathbf{J}_{color}(P)$. 
Let $\maxelt$ denote the unique maximal element of $L$. 
For any $\telt \in L \setminus \{\maxelt\}$, let $k$ be largest in the set $[m]$ such that $\mathcal{C}_{k} \setminus \telt$ is nonempty, and let $v(\telt) := z_{k}$, the periphery vertex in the chain $\mathcal{C}_{k} = \chain^{-1}(k)$. 
In particular, note that $v(\telt) \not\in \telt$ and that $\telt \cup \mathcal{C}_{k}$ is an order ideal from $P$. 
Let $\kappa: L \setminus \{\maxelt\} \longrightarrow \Delta$ be the function defined by $\kappa(\telt) := \color(v(\telt))$. 

Now for the main result of this section. 
When $P$ has the max property, 
\GridPosetColoring\ asserts that the above vertex-coloring function $\kappa: L \setminus \{\maxelt\} \longrightarrow \Delta$ will meet the criteria of \ChainProductVersion. 
The splitting conclusion of \MainTheorem\ is an immediate consequence. 
The proof of \GridPosetColoring\ is somewhat tedious but mainly involves elementary reasoning about finite posets. 

\noindent 
{\bf \GridPosetColoring}\ \  {\sl Let $P$ be a two-color grid poset as above, and assume that $P$ has the max property.  
With} $L := \mathbf{J}_{color}(P)${\sl , define $\kappa: L \setminus \{\maxelt\} \longrightarrow \Delta$ as above.  
Take $\telt \in L \setminus \{\maxelt\}$, and set $\gamma := \kappa(\telt)$. 
Then} $\{\selt \in \comp_{\gamma}(\telt)\, |\, \selt \not= \maxelt \mbox{ and } \kappa(\selt) = \gamma\}$ 
{\sl is a sub-face of $\comp_{\gamma}(\telt)$.}

{\em Proof.} Let $\mathcal{T}_{\gamma}^{(1)}, \ldots , \mathcal{T}_{\gamma}^{(k)}$ be the chains in $P$ of color $\gamma$ as in the paragraphs preceding \ComponentsAreChains. 
Let $j \in \{1, 2, \ldots , k\}$ such that $v(\telt) \in \mathcal{T}_{\gamma}^{(j)}$.  
The set $\mathcal{F} := \{\selt \in \comp_{\gamma}(\telt)\, |\, v(\telt) \in  \selt\}$ is a face of $\comp_{\gamma}(\telt)$. 
To complete the proof it suffices to show that for any $\selt \in \mathcal{F}$ with $\selt \not= \maxelt$ we have $\kappa(\selt) \not= \gamma$ and that for any $\selt \in \comp_{\gamma}(\telt) \setminus \mathcal{F}$ we have $\kappa(\selt) = \gamma$. 

We use a contradiction argument to show that for any $\selt \in \mathcal{F}$ with $\selt \not= \maxelt$ we have $\kappa(\selt) \not= \gamma$. 
To that end, suppose that for some $\selt \in \mathcal{F}$ with $\selt \not= \maxelt$, we have $\kappa(\selt) = \gamma$. 
Let $l \in \{1, 2, \ldots , k\}$ such that $v(\selt) \in \mathcal{T}_{\gamma}^{(l)}$. 
Our contradiction will be to rule out the possibility that $l=j$, $l > j$, or $l < j$. 
If $l=j$, then $v(\selt) = v(\telt)$. 
Since $\selt \in \mathcal{F}$, we have $v(\telt) \in \selt$. 
So, $v(\selt) \in \selt$, which violates the definition of $v(\selt)$. 
Thus, the $l=j$ case is ruled out.  

Next suppose $l > j$. 
We consider two possibilities: $\mathcal{T}_{\gamma}^{(l)}$ is one of the chains $\mathcal{C}_{1}$ or $\mathcal{C}_{2}$, or it is not. 
If so, then the fact that $l > j$ together with the first requirement of the max property means that $\mathcal{C}_{\gamma}^{(l)}$ must be $\mathcal{C}_{2}$, so $\mathcal{T}_{\gamma}^{(j)} = \mathcal{C}_{1}$. 
In particular, the chains $\mathcal{C}_{1}$ and $\mathcal{C}_{2}$ have the same color. 
Now, consecutive chains in a two-color grid poset have opposite colors unless there are no edges between the chains. 
Therefore there are no edges between $\mathcal{C}_{1}$ and $\mathcal{C}_{2}$. 
Then each maximal element of $\mathcal{C}_{1}$ and $\mathcal{C}_{2}$ is a maximal element of $P$. 
But these two maximal elements have the same color, violating the second requirement of the max property. 
So we conclude that $\mathcal{T}_{\gamma}^{(l)}$ cannot be one of the chains $\mathcal{C}_{1}$ or $\mathcal{C}_{2}$. 
Thus, $\chain(v(\selt)) > 2$. 
From the max property, we conclude that the maximal element $v(\selt)$ of $\mathcal{T}_{\gamma}^{(l)}$ is strictly less than the maximal element $z_{2}$ of $\mathcal{C}_{2}$. 
If  $\mathcal{T}_{\gamma}^{(j)} = \mathcal{C}_{1}$, then it must be the case that $z_{2} \in \telt$, since, by definition of $v(\telt)$, every periphery element $z$ with $\chain(z) > \chain(v(\telt))$ is in $\telt$. 
The max property guarantees that $z_{2}$ necessarily has a color different from $\gamma$. 
Since the vertices in $\telt$ which have a color different from $\gamma$ are also in $\selt$, we conclude that $z_{2} \in \selt$ as well. 
Since $\selt$ is an order ideal, the fact that $z_{2} \in \selt$ and $v(\selt) < z_{2}$ forces $v(\selt) \in \selt$, in violation of the fact that $v(\selt)$ is chosen to be an element not in $\selt$. 
So $\mathcal{T}_{\gamma}^{(j)} \not= \mathcal{C}_{1}$. 
Observe that the max property forces $v(\selt)$ to be strictly less than the periphery elements $z_{1}, z_{2}, \ldots z_{(\chain(v(\selt))-1)}$. 
Now the hypothesis $l > j$ means $\chain(v(\telt)) < \chain(v(\selt))$, so the latter list includes $v(\telt)$, hence $v(\selt) < v(\telt)$. 
But since $\selt$ is an order ideal, the fact that $v(\telt) \in \selt$ and $v(\selt) < v(\telt)$ forces $v(\selt) \in \selt$, in violation of the fact that $v(\selt)$ is chosen to be an element not in $\selt$. 
We have therefore ruled out the possibility that $l > j$. 

Now consider the case that $l < j$.  
Since $l < j$, then $\chain(v(\selt)) < \chain(v(\telt))$ and hence $\chain(v(\telt)) \geq 2$. 
As in the previous paragraph, the max property requires that $v(\telt)$ be strictly less than each of the periphery elements $z_{2}, z_{3}, \ldots z_{(\chain(v(\telt))-1)}$. 
We will show that $\chain(v(\selt)) = 1$ by considering the cases $v(\selt) \in \telt$ and $v(\selt) \not\in \telt$; then we will show that $\chain(v(\selt)) = 1$ leads to a contradiction, thereby eliminating the possibility that $l < j$.  
First, assume that $v(\selt) \in \telt$. 
Now, $v(\selt)$ cannot be one of the periphery elements $z_{2}, z_{3}, \ldots z_{(\chain(v(\telt))-1)}$, because this would mean $v(\telt) < v(\selt)$, and the fact that $\telt$ is an order ideal would force $v(\telt) \in \telt$. 
Second, consider the possibility that $v(\selt) \not\in \telt$. 
The fact that $\phi$ from \ComponentsAreChains\ is a poset isomorphism means that $\telt \cup \mathcal{T}_{\gamma}^{(l)}$ is an order ideal. 
This order ideal contains $v(\selt)$ but not $v(\telt)$. 
Therefore $v(\selt)$ cannot be amongst the periphery elements $z_{2}, z_{3}, \ldots z_{(\chain(v(\telt))-1)}$. 
Then $\chain(v(\selt)) = 1$. 
Now $\chain(v(\selt)) = 1$ means in particular that $\selt$ contains every periphery element $z$ with $\chain(z) > 1$. 
So $z_{2} \in \selt$. 
As argued in the previous paragraph, the max property requires that $\color(z_{2}) \not= \color(v(\selt))$, so $\color(z_{2}) \not= \gamma$.  
Since $\selt$ and $\telt$ can only differ in their color $\gamma$ vertices, then $z_{2} \in \telt$. 
But now, $v(\telt) \leq z_{2}$ in $P$ means that we must have $v(\telt) \in \telt$. 
We have therefore ruled out the possibility that $l < j$. 

Thus we have shown that for any $\selt \in \mathcal{F}$ with $\selt \not= \maxelt$, we have $\kappa(\selt) \not= \gamma$. 

Next, we show that for any $\selt \in \comp_{\gamma}(\telt) \setminus \mathcal{F}$, we have $\kappa(\selt) = \gamma$. 
Supposing otherwise, assume $\kappa(\selt) \not= \gamma$ for some $\selt \in \comp_{\gamma}(\telt) \setminus \mathcal{F}$. 
We cannot have $\chain(v(\selt)) < \chain(v(\telt))$: else, the definition of $v(\selt)$ would require that the periphery element $v(\telt)$ be in $\selt$, in violation of the fact that $\selt \not\in \mathcal{F}$. 
So, $\chain(v(\selt)) \geq \chain(v(\telt))$. 
This inequality must be strict, otherwise we get $v(\selt) = v(\telt)$ and $\color(v(\selt)) = \color(v(\telt)) =  \gamma$. 
The fact that $\phi$ from \ComponentsAreChains\ is a poset isomorphism means that $\selt \cup \mathcal{T}_{\gamma}^{(j)}$ is an order ideal. 
This order ideal contains $v(\telt)$ but not $v(\selt)$. 
Therefore $v(\telt)$ cannot be amongst the periphery elements $z_{2}, z_{3}, \ldots z_{(\chain(v(\selt))-1)}$.
Then $\chain(v(\telt)) = 1$. 
Now $\chain(v(\telt)) = 1$ means in particular that $\telt$ contains every periphery element $z$ with $\chain(z) > 1$. 
Then $z_{2} \in \telt$. 
As argued in the previous paragraph, the max property requires that $\color(z_{2}) \not= \color(v(\telt))$, so $\color(z_{2}) \not= \gamma$.  
Since $\selt$ and $\telt$ can only differ in their color $\gamma$ vertices, then $z_{2} \in \selt$. 
But now, $v(\selt) \leq z_{2}$ in $P$ means that we must have $v(\selt) \in \selt$, contradicting the fact that by definition $v(\selt) \not\in \selt$. 
We have therefore ruled out the possibility that $\kappa(\selt) \not= \gamma$, so $\kappa(\selt) = \gamma$.\hfill\QED

{\bf \S \ProofSection\ \ Proof of our main theorem and its corollary.} 

{\em Proof of \MainTheorem.} 
In the notation of the theorem statement, let $\maxelt$ denote the unique maximal element of $L$, so $\lambda = wt(\maxelt)$. 
\ComponentsAreChains\ asserts that each one-color component of $L$ is isomorphic to a product of chains, so that $L$ satisfies one of the hypotheses required for us to be able to invoke \ChainProductVersion. 
Take $\mathcal{S} := \{\maxelt\}$. 
By \GridPosetColoring, the vertex-coloring function $\kappa: L \setminus \mathcal{S} \longrightarrow \{\alpha,\beta\}$ of \S \ColorFunctionSection\ satisfies the vertex-coloring requirements of \ChainProductVersion.  
Now applying \ChainProductVersion, we obtain that $L$ is a splitting distributive lattice for $\chi_{_{\lambda}}$. 
The remaining conclusions of \MainTheorem\ now follow from \SplittingBasicsProposition. 
Computations carried out in \S 5 of \cite{ADLMPPW} for the irreducible rank two root systems show how the expression in \SplittingBasicsProposition\ for the rank-generating function of $L$ becomes the expressions given in the statement of \MainTheorem. 
The formula for the rank generating function in the $\myA_{1} \oplus \myA_{1}$ case is easily derived and is left as a pleasant exercise for the reader.\hfill\QED 

{\em Proof of \MainCorollary.} 
By \StructureMaxResult, we can apply \MainTheorem\ to all $\Phi$-semistandard posets/lattices in order to obtain the conclusions of \MainCorollary.\hfill\QED

\renewcommand{\refname}{\large \bf References}
\renewcommand{\baselinestretch}{1.1}

\end{document}